\documentclass[reqno]{amsart}

\usepackage{amsmath,amsthm,amssymb,amscd}
\usepackage{amsfonts}
\usepackage{rotating}
\usepackage{euscript}
\usepackage{pst-node}
\usepackage{epsfig,verbatim}

\def\be{\begin{equation}}
\def\ee{\end{equation}}

\def\C{{\mathbb C}} 
\def\f{\EuScript}
 
\def\P{{\mathbb P}}

\def\phi{{\varphi}}
\def\v{{\varepsilon}} 
\def\tt{\widetilde}
\def\deg{{\rm deg\,}}

\def\GCD{{\rm GCD }}

\def\bp{\begin{proposition}}
\def\ep{\end{proposition}}

\def\bt{\begin{theorem}}
\def\et{\end{theorem}}
\def\br{\begin{remark}}
\def\er{\end{remark}}
\def\be{\begin{equation}}
\def\bee{\begin{equation*}}
\def\l{\label}

\def\ee{\end{equation}}
\def\eee{\end{equation*}}
\def\bl{\begin{lemma}}
\def\el{\end{lemma}}
\def\bc{\begin{corollary}}
\def\ec{\end{corollary}}
\def\pr{\noindent{\it Proof. }}

\def\bd{\begin{definition}}
\def\ed{\end{definition}}
\def\t{\widetilde}

\newtheorem{theorem}{Theorem}[section]
\newtheorem{lemma}{Lemma}[section]
\newtheorem{definition}{Definition}[section]
\newtheorem{corollary}{Corollary}[section]
\newtheorem{proposition}{Proposition}[section]
\newtheorem{remark}{Remark}[section]

\begin{document}

\author{Fedor Pakovich}

\title{On algebraic curves 
$A(x)-B(y)=0$ of genus zero}

\begin{abstract}
Using a geometric approach involving Riemann surface orbifolds, we provide  lower bounds for the genus of an irreducible algebraic curve  of the form
$\f E_{A,B}:\, 
A(x)-B(y)=0$, where $A, B\in\C(z)$.
We also investigate ``series'' of curves $\f E_{A,B}$ of genus zero, where by a series we mean a family with the ``same'' $A$.  
We show that for a given rational function $A$ 
a sequence of rational functions $B_i$, such that $\deg B_i\rightarrow  \infty$ and all the curves 
$A(x)-B_i(y)=0$ are irreducible and have genus zero, exists 
if and only if  the Galois closure of the field extension  $\C(z)/\C(A)$ has genus zero  or one.
\end{abstract}

\maketitle

\section{Introduction}
\label{intro}

The study
of irreducible algebraic curves of genus zero 
having the form
\be \l{1} \f E_{A,B}:\, A(x)-B(y)=0,\ee
where $A$ and $B$ are complex polynomials, 
has two main motivations. 
On the one hand, such curves 
have special  Diophantine properties.
Indeed, by the Siegel theorem,  
if an irreducible algebraic curve $\f C$ with rational coefficients has infinitely many integer points, then $\f C$ is of genus zero with at most two points at infi\-nity.
More generally, by the Faltings theorem,  if $\f C$ 
has infinitely many rational points, then  $g(\f E_{A,B})\leq 1.$ Therefore, since many interesting Diophantine equations have the form $ A(x)=B(y)$, where $A,B\in \mathbb Q[z],$
the problem of description of curves $\f E_{A,B}$ of genus zero is important for the number theory
(see e.g. \cite{f1}, \cite{bilu},   \cite{kt}).

On the other hand, for polynomials $A$ and $B$ with arbitrary complex coefficients the equality $g(\f E_{A,B})=0$ holds if and only if there exist $C,D\in \C(z)$ satisfying the 
functional equation 
\be \l{2} A\circ C=B\circ D.\ee 
Since equation \eqref{2} describes situations in which a rational function can be decomposed 
into a composition of rational functions in two different ways,  this equation plays a central role in the theory of functional decompositions of rational functions. 
Furthermore, functional equation \eqref{2} where $C$ and $D$ are allowed to be {\it entire} functions reduces to the case where $C,D\in \C(z)$ (see \cite{bn}, \cite{ent}). 
Thus, the problem of description of curves $\f E_{A,B}$ of genus zero naturally appears also in the study of functional equations 
(see e. g. \cite{f1}, \cite{lau}, \cite{ent}, \cite{alg}).

Having in mind possible applications to equation \eqref{2} in rational functions, in this paper we study curves $\f E_{A,B}$ allowing $A$ and $B$ to be arbitrary {\it rational} functions meaning 
by $\f E_{A,B}$ 
the expression obtained by equating to zero the numerator of $A(x)-B(y).$ 
Notice that the curve $\f E_{A,B}$ 
may turn out reducible. In this case its analysis is more complicated and has a different flavor (see e.g. \cite{frii}), 
so below we always will assume that considered curves  $\f E_{A,B}$ are irreducible.

For polynomial $A$ and $B$ the classification of curves ${\f E}_{A,B}$ of genus zero  with one point at infinity 
follows from the so-called ``second Ritt theorem"  (\cite{r1})
about polynomial solutions of \eqref{2}.
Namely, any such a curve has either the form 
\be \l{lau} x^n-y^sR^n(y)=0,\ee where $R$ is an arbitrary polynomial and $\GCD(s,n)=1,$ or the form
\be \l{che} T_n(x)-T_m(y)=0,\ee where $T_n,$ $T_m$ are Chebyshev polynomials and $\GCD(n,m)=1.$
The classification of polynomial  curves $\f E_{A,B}$ of genus zero  with at most two points at infinity
was obtained in the
paper of Bilu and Tichy \cite{bilu}, which continued the line of researches started by Fried (see \cite{f1}, \cite{f2}, \cite{f3}).
In this case, in addition to the above curves we have the following possibilities: 
\be \l{s} x^2-(1-y^2)S^2(y)=0,\ee where $S$ is an arbitrary polynomial,
\be \l{che2} T_{2n}(x)+T_{2m}(y)=0,\ee where  $\GCD(n,m)=1,$ and 
\be \l{last} 
(3x^4-4x^3)-(y^2-1)^3=0.
\ee
Finally, the classification obtained in \cite{bilu} was extended to the case where $A$ and $B$ are allowed to be Laurent polynomials in \cite{lau}. 
In this case, to the list above one has to add the possibility for $R$ in \eqref{lau} to be a Laurent polynomial, and 
the curve 
\be \l{j} \frac{1}{2}\left(y^n+\frac{1}{y^n}\right)-T_m(x)=0,\ee  
where  $\GCD(n,m)=1.$ 
Notice also that an explicit classification of curves \eqref{1} of genus {\it one} with one point at infinity for 
polynomial $A$ and $B$ 
was obtained  by Avanzi and Zannier in \cite{az}. The above results 
essentially exhaust the list of  general results concerning the problem of description of irreducible curves $\f E_{A,B}$ of small genus.

All the curves $\f E_{A,B}$ of genus zero  listed above, except for \eqref{last}, obviously share the following feature: in fact   they are 
``series'' of curves with the ``same'' $A$.
We formalize this observation as follows. Say that a rational function $A$ is a {\it basis of series of curves of genus zero} if there exists a sequence of rational functions $B_i$  such that $\deg B_i\rightarrow  \infty$ and 
all the curves 
$A(x)-B_i(y)=0$ are irreducible  and have genus zero. 
Clearly, a description of all bases of series is an important step in understanding  of the general problem, and the main goal of the paper is to provide such a description in geometric terms.

Recall that for a rational function $A$  
its  normalization $\t A$ is defined as a holomorphic function of the lowest possible degree
between compact Riemann surfaces $\t A:\,\t S_A\rightarrow \C\P^1$  such that $\t A$ is a Galois covering and
$\t A=A\circ H$ for some  holomorphic map $H:\,\t S_A\rightarrow \C\P^1$. From the algebraic point of view, 
the passage from $A$ to $\t A$ 
corresponds to the passage from the field  extension $\C(z)/\C(A)$ to its Galois closure. In these terms our main result about bases of series
is the following statement.

\bt\l{m}
A rational function $A$ is a basis of series of curves of genus zero
if and only if the Galois closure of 
$\C(z)/\C(A)$ has genus zero  or one.
\et

Thus, the set of possible bases of series splits into two classes.  Elements of the first class are  ``compositional left  factor'' of    
well known Galois coverings of $\C\P^1$ by $\C\P^1$ 
calculated for the first time by Klein (\cite{klein}). In particular, up to the change $A\rightarrow \mu_1\circ A\circ \mu_2$, where $\mu_1$ and $\mu_2$ are M\"obius transformations,  besides  
the functions
\be \l{ure} z^n, \ \ \ \ \ T_n, \ \ \ \ \ \frac{1}{2}\left(z^n+\frac{1}{z^n}\right),\ \ \ n\geq1, \ee 
this class contains only a finite number of functions which can be calculated explicitly. 
For instance, the polynomial $3x^4-4x^3$ appearing in  \eqref{last} is an example of such a function, implying that curve \eqref{last} in fact also belongs to a series of curves of genus zero (see Section 5 below). 
Typical representatives of the second class, consisting of rational compositional left factors of  Galois coverings of $\C\P^1$ by a torus, 
are Latt\`es functions
(see e.g. \cite{mil2}), but other possibilities also exist.

The approach of the papers \cite{az}, \cite{bilu}, \cite{lau} to the calculation of $g(\f E_{A,B})$ is based on the formula, given in \cite{f3}, which expresses $g(\f E_{A,B})$ through the
ramifications of $A$ and $B$.  Namely, if $c_1,c_2,\dots c_r$ is a union of critical values of $A$ and $B$, and 
$
f_{i,1},f_{i,2}, ... , f_{i,u_i}
$ (resp. $g_{i,1},g_{i,2}, ... , g_{i,v_i})$ is a collection of local degrees of $A$ (resp. $B$) at the points of 
$A^{-1}(\{c_i\})$ (resp. $B^{-1}(\{c_i\})$), then $g(\f E_{A,B})$ may be calculated as follows:
\be \l{fri}
2-2g(\f E_{A,B}) =
\sum_{i=1}^{r}
\sum_{j_1=1}^{u_{i}} \sum_{j_2=1}^{v_{i}} \GCD(f_{i,j_1}g_{i,j_2})-(r-2)\deg A\,\deg B.
\ee
However, the direct analysis of this formula is quite  difficult
already in the above cases,
and the further progress requires even more cumbersome
considerations. 
In this paper we propose a new approach to the problem
and prove the following general result.

\bt \l{m2} 
Let $A$ be a rational function of degree $n$.
Then for any rational function 
$B$ of degree $m$ such that the curve $\f E_{A,B}$ is irreducible
the inequality 

\be \l{ma} 
g(\f E_{A,B}) >\frac{m-84n+168}{84}
\ee holds, unless  the Galois closure of  $\C(z)/\C(A)$ has genus zero  or one.
\et 

Our approach is based on techniques introduced 
in the recent paper \cite{semi}. This paper studies 
rational solutions of the functional equation \be \l{se} A\circ X=X\circ B\ee  
using Riemann surface orbifolds. For the first time  orbifolds were used in the context of functional equations in the paper 
 \cite{e2} devoted to commuting rational functions. However, in \cite{e2} 
orbi\-folds appear in a dynamical context as a certain characteristic of the Poincar\'e function, while in \cite{semi} an orbifold is attached directly to any rational function.
The approach of \cite{semi} permits to obtain  restrictions on possible ramifications of solutions of \eqref{2} in  terms of the corresponding orbifolds, and to give transparent  proofs of Theorems \ref{m} and  \ref{m2}.

The paper is organized as follows. In the second section 
we recall basic facts about Riemann surface orbifolds and some results from the papers \cite{lau} and \cite{semi}. 
We also express the condition that the Galois closure of  $\C(z)/\C(A)$ has genus zero  or one in terms of orbifolds.
In the third and the fourth sections  we prove Theorem \ref{m2} and Theorem \ref{m} correspondingly. Finally, in the fifth section we consider an  example illustrating 
Theorem \ref{m}.

\section{Fiber products, orbifolds, and Galois coverings}
A pair $\f O=(R,\nu)$ consisting of a Riemann surface $R$ and a ramification function $\nu:R\rightarrow \mathbb N$ which takes the value $\nu(z)=1$ except at isolated points is called a Riemann surface orbifold (see e.g.  \cite{mil}, Appendix E). The Euler characteristic of an orbifold $\f O=(R,\nu)$ is defined by the formula  
\be \l{char}  \chi(\f O)=\chi(R)+\sum_{z\in R}\left(\frac{1}{\nu(z)}-1\right),\ee where $\chi(R)$ is the 
Euler characteristic of $R.$ 
If $R_1$, $R_2$ are Riemann surfaces provided with ramification functions $\nu_1,$ $\nu_2$, and 
$f:\, R_1\rightarrow R_2$ is a holomorphic branched covering map, then $f$
is called  {\it a covering map} $f:\,  \f O_1\rightarrow \f O_2$
{\it between orbifolds}
$\f O_1=(R_1,\nu_1)$ {\it and }$\f O_2=(R_2,\nu_2)$
if for any $z\in R_1$ the equality 
\be \l{us} \nu_{2}(f(z))=\nu_{1}(z)\deg_zf\ee holds, where $\deg_zf$ denotes the local degree of $f$ at the point $z$.
If for any $z\in R_1$ instead of equality \eqref{us} 
a weaker condition 
\be \l{uuss} \nu_{2}(f(z))\mid \nu_{1}(z)\deg_zf\ee
holds,  then $f$
is called {\it a holomorphic map} $f:\,  \f O_1\rightarrow \f O_2$
{\it between orbifolds}.
$\f O_1$ {\it and} $\f O_2.$

A universal covering of an orbifold ${\f O}$
is a covering map between orbifolds \linebreak $\theta_{\f O}:\,
\tt {\f O}\rightarrow \f O$ such that $\tt R$ is simply connected and $\tt \nu(z)\equiv 1.$ 
If $\theta_{\f O}$ is such a map, then 
there exists a group $\Gamma_{\f O}$ of conformal automorphisms of $\tt R$ such that the equality 
$\theta_{\f O}(z_1)=\theta_{\f O}(z_2)$ holds for $z_1,z_2\in \tt R$ if and only if $z_1=\sigma(z_2)$ for some $\sigma\in \Gamma_{\f O}.$ A universal covering exists and 
is unique up to a conformal isomorphism of $\tt R,$
unless $\f O$ is the Riemann sphere with one ramified point, or $\f O$ is the Riemann sphere with two ramified points $z_1,$ $z_2$ such that $\nu(z_1)\neq \nu(z_2)$.
 Furthermore, 
$\tt R=\mathbb D$ if and only if $\chi(\f O)<0,$ $\tt R=\C$ if and only if $\chi(\f O)=0,$ and $\tt R=\C\P^1$ if and only if $\chi(\f O)>0$ 
  (see  \cite{mil}, Appendix E, and \cite{fk}, Section IV.9.12).
Abusing  notation we will use the symbol $\tt {\f O}$ both for the
orbifold and for the  Riemann surface  $\tt R$.

Covering maps between orbifolds lift to isomorphisms between their universal coverings. More generally, the following proposition holds (see \cite{semi}, Proposition 3.1).

\bp \l{poiu} Let $f:\,  \f O_1\rightarrow \f O_2$ be a holomorphic map between orbifolds. Then for any choice of $\theta_{\f O_1}$ and $\theta_{\f O_2}$ there exist 
a holomorphic map \linebreak $F:\, \tt {\f O_1} \rightarrow \tt {\f O_2}$ and 
a homomorphism $\phi:\, \Gamma_{\f O_1}\rightarrow \Gamma_{\f O_2}$ such that the diagram 
\be \l{dia2}
\begin{CD}
\tt {\f O_1} @>F>> \tt {\f O_2}\\
@VV\theta_{\f O_1}V @VV\theta_{\f O_2}V\\ 
\f O_1 @>f >> \f O_2\ 
\end{CD}
\ee
is commutative and 
for any $\sigma\in \Gamma_{\f O_1}$ the equality
\be \l{homm}  F\circ\sigma=\phi(\sigma)\circ F \ee holds.
The map $F$ is defined by $\theta_{\f O_1}$, $\theta_{\f O_2}$, and $f$  
uniquely up to a transformation 
$F\rightarrow g\circ F,$ where $g\in \Gamma_{\f O_2}$. 
In the other direction, for any holomorphic map \linebreak $F:\, \tt {\f O_1} \rightarrow \tt {\f O_2}$  which satisfies \eqref{homm} for some homomorphism $\phi:\, \Gamma_{\f O_1}\rightarrow \Gamma_{\f O_2}$
there exists a uniquely defined  holomorphic map between orbifolds \linebreak $f:\,  \f O_1\rightarrow \f O_2$ such that diagram \eqref{dia2} is commutative.
The holomorphic map $F$ is an isomorphism if and only if $f$ is a covering map between orbifolds. \qed

\ep

If $f:\,  \f O_1\rightarrow \f O_2$ is a covering map between orbifolds with compact support, then  the Riemann-Hurwitz 
formula implies that 
\be \l{rhor} \chi(\f O_1)=d \chi(\f O_2), \ee
where $d=\deg f$. 
For holomorphic maps the following statement is true (see \cite{semi}, Proposition 3.2). 

\bp \l{p22} Let $f:\, \f O_1\rightarrow \f O_2$ be a holomorphic map between orbifolds with compact support.
Then 
\be \l{iioopp} \chi(\f O_1)\leq \chi(\f O_2)\,\deg f \ee and the equality 
holds if and only if $f:\, \f O_1\rightarrow \f O_2$ is a covering map between orbifolds. \qed
\ep

Let $R_1$, $R_2$ be Riemann surfaces, and 
$f:\, R_1\rightarrow R_2$ a holomorphic branched covering map. Assume that $R_2$ is provided with ramification function $\nu_2$. In order to define a ramification function $\nu_1$ on $R_1$ so that $f$ would be a holomorphic map between orbifolds $\f O_1=(R_1,\nu_1)$ and $\f O_2=(R_2,\nu_2)$ 
we must satisfy condition \eqref{uuss}, and it is easy to see that
for any  $z\in R_1$ a minimal possible value for $\nu_1(z)$ is defined by 
the equality 
\be \l{rys} \nu_{2}(f(z))=\nu_{1}(z)\GCD(\deg_zf, \nu_{2}(f(z)).\ee 
In case if \eqref{rys} is satisfied for  any $z\in R_1$ we 
say that $f$ is  {\it a minimal holomorphic  map 
between orbifolds} 
$\f O_1=(R_1,\nu_1)$ {\it and} $\f O_2=(R_2,\nu_2)$. Notice that any covering map obviously is a minimal holomorphic  map.

With any holomorphic function $f:\, R_1\rightarrow R_2$ between compact Riemann surfaces 
one can associate in a natural way two orbifolds $\f O_1^f=(R_1,\nu_1^f)$ and 
$\f O_2^f=(R_2,\nu_2^f)$, setting $\nu_2^f(z)$  
equal to the least common multiple of local degrees of $f$ at the points 
of the preimage $f^{-1}\{z\}$, and $$\nu_1^f(z)=\nu_2^f(f(z))/\deg_zf.$$ By construction, 
$f$ is a covering map between orbifolds $f:\, \f O_1^f\rightarrow \f O_2^f.$ 
Furthermore, since the composition $f\circ \theta_{\f O_1^f}: \t{\f O_1^f}\rightarrow \f O_2^f$ is a covering map  between orbifolds, it follows from the uniqueness of the universal covering that 
\be \l{ravv}  \theta_{\f O_2^f}=f\circ \theta_{\f O_1^f}.\ee

For rational functions $A$ and $B$ irreducible components of  ${\f E}_{A,B}$ correspond to irreducible components of the 
fiber product of $A$ and $B.$ In particular, if ${\f E}_{A,B}$  is an irreducible
curve and $\t{\f E}_{A,B}$ is its desingularization, then there exist holomorphic functions $p,q:\t{\f E}_{A,B}\rightarrow \C\P^1$ such that 
\be \l{xui} 
A\circ p=B\circ q,
\ee
and \be \l{xer} \deg A=\deg q,\ \ \ \deg B=\deg p\ee (see \cite{lau}, Theorem 2.2 and Proposition 2.4).
Furthermore, the functions $A,B,$ $p,q$ possess ``good'' properties with respect to the associated orbifolds defined above. 
Namely, 
the following statement  holds (see \cite{semi}, Theorem 4.2 and Lemma 2.1).

\bt \l{t1} Let $A$, $B$ be rational functions such that the curve ${\f E}_{A,B}$ is irreducible, and 
$p,q:\t{\f E}_{A,B}\rightarrow \C\P^1$ holomorphic functions such that 
equalities \eqref{xui} and \eqref{xer} hold.  Then 
the commutative diagram 
$$ 
\begin{CD}
\f O_1^q @>p>> \f O_1^A\\
@VV q V @VV A V\\ 
\f O_2^q @>B >> \f O_2^A\ 
\end{CD}
$$
consists of minimal holomorphic  maps between orbifolds. \qed
\et

Of course, vertical arrows in the above diagram are covering maps and hence  minimal holomorphic  maps  simply by definition. The meaning of the theorem is that 
the branching of $q$ and $A$ to a certain extent defines 
the branching of $p$ and $B$. 
For example, Theorem \ref{t1} 
applied  to functional equation \eqref{se}
where $A,X,B$ are rational functions such 
that $\f E_{A,X}$ is irreducible, implies that  $\chi(\f O_2^X)\geq 0$ 
(see   \cite{semi}). 

For a rational function $A$ the condition $\chi(\f O_2^A)\geq 0$ is very restrictive, and  is  equivalent to the condition 
that the normalization of $A$  has genus at most one. 

\bl \l{ml} 
Let $A$ be a  rational function. Then $g(\t S_A)=0$ if and only if \linebreak $\chi(\f O_2^A)> 0$, and  $g(\t S_A)=1$ if and only if  $\chi(\f O_2^A)= 0$. 
\el
\pr 
Let $f: S\rightarrow \C\P^1$ be an arbitrary Galois covering of $\C\P^1$. Then $f$ is a quotient map 
$f:S\rightarrow S/\Gamma$ for some subgroup $\Gamma$ of $Aut(S)$, and for any branch point $z_i$, $1\leq i \leq r,$  of $f$ there exists a number $d_i$ such that $f^{-1}\{z_i\}$ consists of  $\vert G\vert /d_i$ points, at each of which the multiplicity of $f$ equals $d_i$. Applying
the Riemann-Hurwitz formula, we see that 
$$2g(S)-2=-2\vert \Gamma\vert +\sum_{i=1}^r\frac{\vert \Gamma\vert}{d_i}\left(d_i-1\right),$$
implying that
\be \l{gc} \chi(\f O_2^f)=2+\sum_{i=1}^r\left(\frac{1}{d_i}-1\right)=\frac{2-2g(S)}{\vert \Gamma\vert }.\ee 
Thus, if $f: S\rightarrow \C\P^1$ is a Galois covering, then $g( S)=0$ if and only if $\chi(\f O_2^f)>0$, while  $g( S)=1$ if and only if $\chi(\f O_2^f)=0$.

Let now $A:\, \C\P^1\rightarrow \C\P^1$ be an arbitrary rational function. Since the normalization
$\t A:\, \t S_A\rightarrow \C\P^1$ of $A$ can be described as any  irreducible
component of the $m$-fold  fiber product of $A$ distinct from the diagonal components 
where two or more coordinates are equal 
(see \cite{f}, $\S$I.G), it follows from 
the construction of the fiber product (see e. g. \cite{lau}, Section 2 and 3) that 
\be \l{vtgh} \f O_2^A=\f O_2^{\t A}.\ee  Thus, $g(\t S_A)=0$ if and only if  $\chi(\f O_2^A)> 0$, and  $g(\t S_A)=1$ if and only if \linebreak $\chi(\f O_2^A)= 0$. 
\qed

\vskip 0.2cm

If $\f O=(\C\P^1, \nu)$ is an orbifold such that 
$\chi(\f O)= 0$, then \eqref{char} implies that  the collection of ramification indices of $\f O$ is either $(2,2,2,2)$, or one of the following triples $(3,3,3)$, $(2,4,4)$, $(2,3,6)$. 
For all such orbifolds $\tt{\f O}=\C$. Furthermore, the group $\Gamma_{\f O}$
is generated by translations of $\C$ by elements of some lattice $L\subset \C$ of rank two and the transformation $z\rightarrow  \v z,$ where $\v$ is $n$th root of unity with $n$ equal to 2,3,4, or 6, such that  $\v L=L.$ 
For the collection of ramification indices $(2,2,2,2)$ the complex structure of $\C/L$ may be arbitrary and the function  $\theta_{\f O}$ is the corresponding Weierstrass function $\wp(z).$ On the other hand, for the collections $(2,4,4)$, $(2,3,6),$ $(3,3,3)$ this structure is   
rigid and arises from the tiling of $\C$ by squares, equilateral triangles, or alternately colored equilateral triangles, respectively. Accordingly, the functions $\theta_{\f O}$ may be written in terms of the corresponding
Weierstrass functions as $\wp^2(z),$ $\wp^{\prime 2}(z)$, and $\wp^{\prime }(z)$   (see  \cite{mil2} and \cite{fk}, Section IV.9.12).

Similarly, if $\chi(\f O)> 0$,  then the collection of ramification indices of $\f O$ is either $(n,n)$ for some $n\geq 2,$ or $(2,2,n)$ for some $n\geq 2,$ or one of the following triples $(2,3,3)$, $(2,3,4),$ $(2,3,5)$. 
In fact, formula \eqref{char}  also allows $\f O$ to be a non-ramified 
sphere or one of two orbifolds without universal covering.
However, if $\f O=\f O_2^A$ for some rational function $A$, then these cases are impossible since 
for any rational function $A$ both orbifolds 
$\f O_1^A$, $\f O_2^A$ have a universal covering (see \cite{semi}, Lemma 4.2), and  $\f O_2^A$ cannot be non-ramified. 
Further, $\tt{\f O}=\C\P^1$, and  
the group $\Gamma_{\f O}$ is a finite subgroup of the automorphism group of 
$\C\P^1$. Namely, to orbifolds with the collections of ramification indices  $(n,n)$, $(2,2,n)$, $(2,3,3)$, $(2,3,4),$ and $(2,3,5)$ correspond the groups $C_n,$ $D_{2n},$ $A_4,$ $S_4,$ and $A_5$.
The corresponding functions  $\theta_{\f O}$ are Galois coverings of $\C\P^1$ by $\C\P^1$ and have degrees $n$, $2n$, 12, 24, and 60 (see \cite{klein}).

\section{Proof of Theorem \ref{m2}}
First of all, observe that if $f:\, R\rightarrow \C\P^1$ is a holomorphic function of degree $n$ on a Riemann  surface $R$ of genus $g$, 
then 
\be  \l{cha} \chi(\f O_2^f)> 4-2g-2n.\ee
Indeed, it follows from the definition that 
$$\chi(\f O_2^f)> 2- c(f),$$ where $c(f)$ denotes the number of branch points of $f.$ On the other 
hand, since the number $c(f)$ is less than or equal to the number of points $z\in R$ where $\deg_z f>1$,  the Riemann-Hurwitz formula 
$$\chi(R)=\chi(\C\P^1)n-\sum_{z\in R}(\deg_zf-1) $$ 
implies 
that $$c(f)\leq \chi(\C\P^1)n-\chi(R).$$ 
Thus, 
$$\chi(\f O_2^f)>2 +\chi(R)-\chi(\C\P^1)n ,$$ implying \eqref{cha}.

Let now $p,q:\,\t{\f E}_{A,B}\rightarrow \C\P^1$ be holomorphic functions  such that \eqref{xui} and \eqref{xer} hold.
Since $B:\, \f O_2^q\rightarrow \f O_2^A$ is a minimal holomorphic map between orbifolds by Theorem \ref{t1},  it follows 
from Proposition \ref{p22} that 
\be \l{bur} \chi(\f O_2^q)\leq m \chi(\f O_2^A).\ee 
On the other hand, \eqref{char} implies that  if  $\chi(\f O)< 0,$ then in fact
\be \l{42} \chi(\f O)\leq -\frac{1}{42}\ee  
(where the equality is attained for the collection of ramification indices $(2,3,7)$).
Therefore, if $\chi(\f O_2^A)<0,$ then \eqref{42} and \eqref{cha} imply  
the inequality 
$$4-2g-2n<-\frac{m}{42}$$ which in  turn implies \eqref{ma}. \qed

\section{Proof of Theorem \ref{m}}
It follows from Theorem \ref{m2} and Lemma \ref{ml} that we only need to show that if  $\chi(\f O_2^A)\geq 0 $, then $A$ is a basis of series. 
Assume first that $\chi(\f O_2^A)= 0$. Then the universal covering of $\f O_2^A$ is $\C$, and 
the group $\Gamma_{\f O_2^A}$
is generated by translations of $\C$ by elements of some lattice $L=<\omega_1,\omega_2>$ and the transformation $z\rightarrow  \v z,$ where $\v$ is an $n$th root of unity with $n$ equal to 2,3,4, or 6, such that  $\v L=L.$ 
This implies that 
for any integer $m\geq 2$ the map $F:z\rightarrow mz$ satisfies condition \eqref{homm}
for the homomorphism $\phi: \Gamma_{\f O_2^A}\rightarrow \Gamma_{\f O_2^A}$ defined on the generators of $\Gamma_{\f O_2^A}$ by the equalities 
\be \l{ho}\phi(z+\omega_1)=z+m\omega_1, \ \ \ \phi(z+\omega_1)=z+m\omega_1,\ \ \  \phi(\v z)=\v z.\ee 
Therefore, by Proposition \ref{poiu},  there exists a rational functions $R_{m}$
such that 
$$\theta_{\f O_2^A}(mz)=R_{m}\circ \theta_{\f O_2^A},$$ and it is easy to see that $\deg R_m=m^2$. 
Furthermore, it follows from \eqref{rhor} that  $\chi(\f O_1^A)= 0$, implying that the group
$\Gamma_{\f O_1^A}$ is generated  by translations by elements of some sublattice $\t L$ of $L$  and the transformation $z\rightarrow \v^l z$ for some $l\geq 1.$ Thus, homomorphism \eqref{ho} satisfies the condition 
\be \l{hh} \phi(\Gamma_{\f O_1^A})\subseteq \Gamma_{\f O_1^A}\,,\ee implying that   
there exists a rational function  $S_{m}$
of degree $m^2$ such that 
$$ \theta_{\f O_1^A}(mz)=S_{m}\circ \theta_{\f O_1^A}.$$ 

Since \be \l{raww} \theta_{\f O_2^A}=A\circ \theta_{\f O_1^A}\,,\ee it follows now from the equalities 
$$\theta_{\f O_2^A}(mz)=R_m \circ \theta_{\f O_2^A}=R_m \circ A\circ\theta_{\f O_1^A} $$
and 
$$\theta_{\f O_2^A}(mz)=A\circ \theta_{\f O_1^A}(mz)=A\circ S_m \circ \theta_{\f O_1^A},$$
 that 
$$A\circ S_m=R_m\circ A.$$ Thus, whenever the curve 
$A(x)-R_m(y)=0$ is irreducible, it has genus zero. Since ${\f E}_{A,B}$ is irreducible whenever the degrees of $A$ and $B$ are coprime (see e. g. \cite{lau}, Proposition 3.1),  
taking any sequence $m_i\to \infty$ whose elements are coprime with $\deg A$, we obtain a sequence  $A(x)-R_{m_i}(y)=0$ of irreducible curves of genus zero.

In the case  $\chi(\f O_2^A)> 0$ the proof is similar with appropriate modifications. First observe that in order to prove the theorem it is enough to show that  for any $A$ with $\chi(\f O_2^A)> 0$ 
there exists a {\it single} pair of rational functions $S$ and $R$ such that 
\be \l{pi} A\circ S=R\circ A\ee and \be \l{xorr} \GCD(\deg R,\deg A)=1.\ee Indeed, 
\eqref{pi} implies that $$ A\circ S^{\circ l}=R^{\circ l}\circ A.$$ Therefore, since equality \eqref{xorr} implies the equality  $\GCD(\deg R^{\circ l},\deg A)=1,$ the sequence  $A(x)-R^{\circ l}(y)=0$ consists of irreducible curves of genus zero.
Further, since by Lemma \ref{ml} the group $\Gamma_{\f O_2^A}$ belongs to the list $C_n,$ $D_{2n},$ $A_4,$ $S_4,$ $A_5$, 
in order to show the existence of such a  pairs 
for any $A$ with $\chi(\f O_2^A)> 0$ it is enough to show that for any group $\Gamma$ from the above list 
there exists a rational function $F$ of degree corpime with $\vert \Gamma \vert$ which is
$\Gamma$-equivariant, that is   satisfies the equality
\be \l{homm+}  F\circ\sigma=\sigma\circ F \ee   for any $\sigma\in \Gamma$. Indeed, condition \eqref{homm+} means  that the corresponding homomorphism in \eqref{homm} satisfies  
$\phi(\sigma)=\sigma$ for any $\sigma \in \Gamma,$ implying that $\phi(\t \Gamma)=  \t \Gamma$ for any subgroup  $\t \Gamma$ of $\Gamma,$ and we conclude 
as above that 
\be \l{w} \theta_{\f O_2^A}\circ F=R\circ \theta_{\f O_2^A}, \ \ \ \theta_{\f O_1^A}\circ F=S\circ \theta_{\f O_1^A}\ee
for some rational functions $S$ and $R$  such that \eqref{pi} holds.
Moreover, since \linebreak $\deg \theta_{\f O_2^A}=\vert \Gamma_{\f O_2^A}\vert $ and $\deg R=\deg F$,  it follows from \eqref{raww} that equality \eqref{xorr} holds.

If $\Gamma_{\f O_2^A}=C_n$, then up to the change $A\rightarrow \mu_1\circ A\circ \mu_2$, where 
$\mu_1,$ $\mu_2$ are M\"obius transformations, $A=z^n$, and hence \eqref{lau} already provides a necessary series of irreducible curves of genus zero.   
Similarly,  if $\Gamma_{\f O_2^A}=D_n$, then without loss of generality we may assume that either $A=T_n$ or $$A=\frac{1}{2}\left(z^n+\frac{1}{z^n}\right)$$ 
(see e.g. Appendix of \cite{mp2}), and hence the statement of the lemma  follows from equalities \eqref{che} and \eqref{j}. 
Finally, since $A_4\subset S_4 \subset A_5$, in order to finish the proof it is enough to find a single $A_5$-equivariant function whose order is coprime with 60, and as such a function we can take for example the function  
\be \l{f} F=\frac{z^{11}+66z^6-11z}{-11z^{10}-66z^5+1}\ee of degree 11, constructed in the paper \cite{dm}.

\section{Example}
Consider the rational function $A=3z^4-4z^3$ appearing in \eqref{last}. The critical values of this function are $0,-1,\infty$. The preimage of $0$ consists of a critical point $0$, whose multiplicity is 3, and the point $4/3.$ The preimage of $-1$ consists of a critical point $1$, whose multiplicity is 2, and the points $-\frac{1}{3}\pm i\frac{\sqrt{3}}{2}.$
 Finally,  the preimage of $\infty$ consists of a single point $\infty$, whose multiplicity is 4. Thus, 
$$\nu_2^A(-1)=2,\ \ \ \nu_2^A(0)=3,\ \ \ \nu_2^A(\infty)=4,$$ and the value of $\nu_2^A$ at any other point equals 1. Correspondingly,  
$$\nu_1^A\left(-\frac{1}{3}+ i\frac{\sqrt{3}}{2}\right)=\nu_1^A\left(-\frac{1}{3}- i\frac{\sqrt{3}}{2}\right)=2, \ \ \ \nu_1^A\left(\frac{4}{3}\right)=3.$$
Finally,  
$$\chi(\f O_2^A)=\frac{1}{12}, \ \ \  \chi(\f O_1^A)=\frac{1}{3},$$ and $\Gamma_{\f O_2^A}=S_4$.

Fix the generators of $S_4$ as 
$$z\rightarrow iz, \ \ \ \   z\rightarrow\frac{z+i}{z-i}.$$ Then $$\theta_{\f O_2^A}=-\frac{(z^8+14z^4+1)^3}{108z^4(z^4-1)^4}.$$  The critical values of $\theta_{\f O_2^A}$  normalized in such a way are 
 $0,-1,\infty$, 
and  $\theta_{\f O_2^A}=A\circ \theta_{\f O_1^A}$, where 
$$\theta_{\f O_1^A}=
\frac{\left(\frac{1}{6}(1+i){z}^{2}-\frac{i}{3}z+ \frac{1}{6}(1-i) \right)\left( {z}^{4}+2\,{z}^{3}+2\,{z}^{2}-2\,z+1 \right)}{\left( {z}^{2}+1 \right)  \left( z+1 \right)  \left( z-1 \right) z}\,.
$$

As an $S_4$-invariant function of degree corpime with $\deg A=4$ we can take function  \eqref{f}. However, we also can take the function of lesser degree 
$$F=\frac{-z^5+5z}{5z^4-1}$$
obtained from  the invariant form 
$$x^5y-xy^5$$ 
by the method of \cite{dm}. For such $F$ the functions $R$ and $S$ from equalities \eqref{w} are
$$R=\frac {{z}^{2} \left( {z}^{3}-240\,{z}^{2}+19200\,z-512000 \right) }{1048576+625\,{z}^{4}+16000\,{z}^{3}+153600\,{z}^{2}\\
\mbox{}+655360\,z}$$
and 
$$S=-\frac {{z}^{2} \left( 3\,{z}^{3}-10\,{z}^{2}+20\,z-40 \right) }{32-20\,{z}^{3}+15\,{z}^{4}}\,.$$ 
Thus, we obtain a family 
of irreducible curves of genus zero 

$$(3x^4-4x^3)-\left(\frac {{y}^{2} \left( {y}^{3}-240\,{y}^{2}+19200\,y-512000 \right) }{1048576+625\,{y}^{4}+16000\,{y}^{3}+153600\,{y}^{2}\\
\mbox{}+655360\,y}\right)^{\circ k}=0,$$
having the parametrizations 
$$x=\left(-\frac {{t}^{2} \left( 3\,{t}^{3}-10\,{t}^{2}+20\,t-40 \right) }{32-20\,{t}^{3}+15\,{t}^{4}}\right)^{\circ k}\hskip -0.3cm, \hskip 2cm y=3t^4-4t^3.$$

\vskip 0.2cm
\noindent{\bf Acknowledgments}. 
 The author is grateful to 
M. Fried  for discussions, and to the Max-Planck-Institut fuer Mathematik for the hospitality and the support.



\end{document}